\documentclass{amsart}
\usepackage{amscd,amsmath}

\newcommand{\D}{\mathcal{D}}
\newcommand{\V}{\mathcal{V}}
\newcommand{\Hh}{\mathcal{H}}

\newcommand{\ov}{\overline}

\newcommand{\Ker}{\mathop\mathrm{Ker\,}}

\numberwithin{equation}{section}

\newtheorem{te}{Theorem}[section]

\newtheorem{pr}{Proposition}[section]
\newtheorem{co}{Corollary}[section]
\newtheorem{lm}{Lemma}[section]

\theoremstyle{definition}
\newtheorem{de}{Definition}[section]
\newtheorem{re}{Remark}[section]
\newtheorem{ex}{Example}[section]

\begin{document}

\title{Holomorphicity and Walczak formula on Sasakian manifolds}

\author{Vasile Br\^{\i}nz\u anescu}
\address{Institute of Mathematics "Simion Stoilow", Romanian Academy -
P.O. Box 1-764, RO-70700 Bucharest, Romania.}
\email{vasile.brinzanescu@imar.ro}

\author{Radu Slobodeanu}
\address{Faculty of Physics, Bucharest University - 405 Atomistilor str., CP Mg-11, RO-76900
Bucharest, Romania.}
\email{slobyr@yahoo.com}

\begin{abstract}
Walczak formula is a very nice tool for understanding the geome\-try of a Riemannian manifold 
equipped with two orthogonal complementary distributions. Svensson \cite {sve} has shown that this 
formula simplifies to a Bochner type formula when we are dealing with K\"ahler manifolds and 
holomorphic (integrable) distributions. We show in this paper that such results have a counterpart in Sasakian geometry. 
To this end, we build on a theory of (contact) holomorphicity on almost contact metric manifolds. Some other applications 
for (pseudo) harmonic morphisms
 on Sasaki manifolds are outlined.
\end{abstract}

\maketitle



\section{Introduction}
Throughout this paper $M$, $N$ etc. will be  connected, $\mathcal{C}^\infty$ manifolds. All geometric objects considered will 
also be smooth.

The analogue of an almost Hermitian structure on odd dimensional spaces is the {\it almost
contact metric structure.} We recall the necessary definitions, cf. \cite{ble}:

\begin{de}
An almost contact structure on a $2n+1$-dimensional manifold $M$ is a triple $(\phi, \xi, \eta)$ where $\phi$ is a  (1, 1) tensor field,
$\xi$ is a vector field  and $\eta$ is a  1-form $\eta$ satisfying  the following relations: 
$$
\phi^{2}=-I+\eta \otimes \xi, \quad \eta(\xi)=1.
$$
A manifold $M$ together with an almost contact structure is called an \emph{almost contact manifold}. $\xi$ is called the 
\emph{characteristic} vector field.

An \emph{almost contact metric structure} $(\phi, \xi, \eta, g)$ is an almost contact
structure together with a compatible metric (which always exists), that is a metric $g$
satisfying: 
$$
g(\phi X, \phi Y) = g(X, Y)-\eta(X)\eta(Y).
$$

If, in addition, $\eta$ is a {\it contact form}
(i.e. $\eta \wedge (d\eta)^{n} \neq 0$) and $g$ is an \emph{associated metric} (i.e. 
$d\eta(X, Y)=g(X, \phi Y)$), then our structure is a \emph{contact metric structure}. 
In this case $\xi$ coincides with the {\it Reeb field} of the {\it contact form} $\eta$.

A contact metric structure whose (1,1)-tensor $\phi$ is normal:
\begin{equation} \label{normalitate}
[\phi, \phi](X, Y) + 2 d \eta(X, Y)\xi = 0
\end{equation}
is called \emph{Sasakian}.

\end{de}

Sasakian structures are the analogue of K\"ahler structures on odd-dimensional manifolds.
 The Sasakian condition is equivalent to the integrability 
of the corresponding almost complex structure on the riemannian cone
 over $M$, cf. \emph{e.g.} \cite{gal}.

The normality equation \eqref{normalitate} is equivalent to the following one:
\begin{equation} \label{Sasaki}
(\nabla_{X}\phi)Y = g(X,Y)\xi-\eta(Y)X,
\end{equation}
which makes transparent the 
analogy with the K\"ahler case: indeed, it is enough to take in both members 
of \eqref{Sasaki} the component tangent to the {\it contact distribution} $\D =\Ker \eta$, 
for $X, Y \in \Gamma(\D)$
and then we obtain a parallelism-type condition
for $\phi$. This is in fact the \emph{transversally K\"ahler} condition.

An almost contact structure has a natural \emph{transversal holomorphic structure}, 
transversality being here understood with respect to the foliation defined by the 
characteristic field. In the language of $G$-structure, this is a $H^{1,n}$-structure, 
cf. \cite{tano}. 

The paper is organized as follows: in \S 2 we  study invariant (to the action of $\phi$) 
distributions on almost contact manifolds. In \S 3 we study a notion of holomorphic 
distribution (in particular, holomorphic vector field), which is automatically 
$\phi$-invariant. We show how is this notion related to holomorphicity on the cone. 
\S 4 is devoted to holomorphicity on normal almost contact manifolds, especially on 
Sasakian manifolds. Finally, in \S 5 we apply our theory of holomorphicity to derive 
results in Riemannian geometry: 
applications of the Walczak formula and properties of some particular harmonic morphisms.

\section{Invariant distributions on almost contact metric manifolds}

In analogy with the definition of a complex distribution on an almost hermitian manifold we give:
\begin{de} Let  $(M, \phi, \xi, \eta, \emph{g})$ be an almost contact metric manifold. 
A distribution $\mathcal{V}$ on $M$ is called {\it invariant} if 
$\phi(\mathcal{V}) \subseteq \mathcal{V}$.
\end{de}

\begin{re}
1. $\D :=\Ker \eta$ is an invariant distribution.

2. On an almost contact metric manifold, a distribution $\mathcal{V}$ is invariant 
if and only if its orthogonal complementary distribution $\mathcal{H}$ is also invariant.
\end{re}

\noindent
The proof follows from the anti-symmetry of $\phi$. Let $X \in \Gamma(\Hh )$, $V \in \Gamma(\V )$; we have:
\begin{equation*}
\begin{split}
{g}(\phi{X},{V})&={g}(\phi^{2}{X}, \phi{V})+
\eta(\phi{X})\eta({V})={g}(-{X}+\eta({X})\xi,\phi{V})\\
&=-{g}({X},\phi{V})+\eta({X})\eta(\phi{V})=
-{g}({X},\phi{V}).
\end{split}
\end{equation*}
By hypothesis, $\phi V \in \Gamma(\V )$, so the last term is zero, which implies that 
$\phi X$ is orthogonal to $V$, for every $V \in \Gamma(\V )$. 
This means $\phi X \in \Gamma(\mathcal{H})$.

\bigskip

Note that, unlike in the Hermitian case, an invariant distribution can be  even 
or odd-dimension as well. In particular, the dimensions of two complementary invariant distributions 
on $M^{2n+1}$ cannot have the same parity.

The position of the characteristic field $\xi$ with respect to an invariant distribution is subject to some restrictions:
\begin{lm}
On an almost contact metric manifold with an invariant distribution $\mathcal{V}$,
 the vector field  $\xi$ must be in $\Gamma(\mathcal{V})$ or in $\Gamma(\mathcal{H})$,
 where $\mathcal{H}=\mathcal{V}^{\perp}$.

Moreover, if $\xi\in\Gamma(\mathcal{V})$, then $\mathcal{H}\subseteq\mathcal{D}$.
\end{lm}
\begin{proof}
Let $\xi^\mathcal{H}$, $\xi^\mathcal{V}$ denote the $\mathcal{H}$, resp. $\mathcal{V}$ component of $\xi$ 
(the exponent $\mathcal{V}$ or $\mathcal{H}$
will always indicate the orthogonal projections onto these distributions).  Then $0=\phi\xi$ together with 
the invariance of $\mathcal{H}$ and $\mathcal{V}$ imply  
$\phi\xi^{\mathcal{H}}=0$, $\phi\xi^{\mathcal{V}}=0$.  
But $\Ker \phi$ is one-dimensional and therefore, if $\xi^{\mathcal{H}}$ and $\xi^{\mathcal{V}}$ were both non-zero, they 
would be collinear, contradiction.

The second statement follows from  $\eta(X)=g(X, \xi)=0$, 
for all ${X}\in\Gamma(\mathcal{H})$.
\end{proof}

On the other hand, the characteristic vector field $\xi$ is tangent to any invariant submanifold of a
 contact metric manifold (cf. \cite[p. 122]{ble}), so one expects the same phenomenon to occur
 for (integrable) invariant distributions. We have indeed:

\begin{pr}
On a contact metric manifold $M^{2n+1}$ endowed with an invariant distribution $\mathcal{V}$
 any of the following conditions implies $\xi\in\Gamma(\mathcal{V})$:

$(i)$ $\mathrm{dim}(\mathcal{V})=2k+1$, \ $k\leq{n}$;

$(ii)$ $\mathcal{V}$ is integrable.

In particular, an integrable invariant distribution must be odd-dimensional. 
\end{pr}

\begin{proof}
$(i)$ By Lemma 2.1, it is enough to prove that $\xi$ is not in $\Gamma(\mathcal{H})$.

 If $\xi\in\Gamma(\mathcal{H})$, then $\mathcal{H}$ admits (local) frames of the
 type $\{\xi, X_i, \phi(X_i) \}$, so it is odd-dimensional, like
 $\mathcal{V}$, contradiction.

$(ii)$ Suppose that $\xi\in\Gamma(\mathcal{H})$. Then, from Lemma 2.1,
 $\mathcal{V} \subseteq \mathcal{D}$, where $\mathcal{D}$ is the contact distribution.
 So, for any $V, W\in\Gamma(\mathcal{V})$, we have
$$
{g}({V},\phi{W})=d\eta({V},{W})=\frac{1}{2}\left[ {V}\eta({W})-{W}\eta({V})-\eta([{V},{W}])\right]=-\frac{1}{2}\eta([{V},{W}])=0,
$$
the last equality being a consequence of the integrability of $\mathcal{V}$.
 We conclude that $\phi{W}$ is orthogonal to $\mathcal{V}$, contradiction.
Hence $\xi$ cannot be in $\Gamma(\mathcal{H})$.
As Lemma 2.1 shows also that $\xi$ can be neither a "mixed" sum, the proof is complete.

(Note that we have not used all the contact structure information, but only that $g$ is an associated metric.) 
\end{proof}

\begin{ex}
On $\mathbb{R}^{2n+1}$ with the standard contact metric structure,
 the distribution $\mathcal{V}_{k}$ ($k \leq n$) locally spanned by 
$$
X_{i}=2\frac{\partial}{\partial y^{i}},\quad 
X_{n+i}=2\left(\frac{\partial}{\partial x^{i}}+y^{i}\frac{\partial}{\partial z}\right)
\quad \text{and possibly} \ \xi \quad (i=\ov{1,k})
$$ 
is an invariant distribution 
of dimension $2k$, respectively $2k+1$ if it contains $\xi$.

\end{ex}

For further use we next prove a relation between the Lie derivative and the covariant derivative of the tensor
 $\phi$, similar to the relation (3.1) in \cite {sve}. The following relation is easily derived:
$$
g(\nabla_{\phi Z}X, V)=g(X, (\mathcal{L}_{V}\phi-\nabla_{V}\phi)Z)-g(X, \phi \nabla_{Z} V).
$$

Using here the anti-symmetry of $\phi$, the fact that $\nabla$ is
a metric connection and also $g(\phi X, V)=0$ (because $\mathcal{H}$ is an invariant 
distribution), we prove:

\begin{pr}
Let $(M, \phi, \xi, \eta, g)$ be an almost contact metric manifold and $\mathcal{V}$
 an invariant distribution with orthogonal complement $\mathcal{H}$. For any section $X$ 
of $\mathcal{H}$ and any vector field $V$ tangent to $\mathcal{V}$, we have:
\begin{equation} \label{svenson}
g(\nabla_{\phi Z}X+\nabla_{Z}\phi X, V)=g(X, (\mathcal{L}_{V}\phi-\nabla_{V}\phi)Z), \quad 
\forall Z \in \Gamma(TM).
\end{equation}
\end{pr}

\medskip
We recall here
that {\it the second fundamental form} $B^{\mathcal{V}}$ and {\it the integrability tensor}
$I^{\mathcal{V}}$, of $\mathcal{V}$, are defined by:
$$
B^{\mathcal{V}}(V, W)=\frac{1}{2}\left(\nabla_{V}W+\nabla_{W}V\right)^{\mathcal{H}},\quad 
I^{\mathcal{V}}(V, W)=[V, W]^{\mathcal{H}}, \quad V, W \in \Gamma(\mathcal{V}).
$$

As for the distribution $\D$, which is invariant, we have:
\begin{re}
On a contact metric manifold,
$$
B^{\D}(X, \phi Y)=B^{\D}(\phi X, Y), \quad \forall X, Y \in \Gamma(\D).
$$
In particular, $\D$ is a minimal distribution. If, in addition, the manifold is
K-contact, then $\D$ is a totally geodesic distribution. 
\end{re}
\begin{proof}
A result of Olszak, \cite {ol}, states that on a contact metric manifold, we have:
\begin{equation}
\left(\nabla_{X}\phi \right)Y+\left(\nabla_{\phi X}\phi \right)\phi Y=2g(X, Y)\xi-
\eta(Y)(X+hX+\eta(X)\xi).
\end{equation}
In particular, if $X, Y \in \Gamma(\D)$, the above relation becomes:
$$
\nabla_{X}\phi Y-\phi \nabla_{X}Y-\nabla_{\phi X}Y-\phi \nabla_{\phi X}\phi Y=2g(X, Y)\xi.
$$
Interchanging $X$ and $Y$, we obtain a similar relation which, subtracted 
from the above one gives:
$$
\nabla_{X}\phi Y+\nabla_{\phi Y}X-(\nabla_{\phi X}Y+\nabla_{Y}\phi X)=
\phi([X, Y]+[\phi X,\phi Y]).
$$
Taking only the component collinear with $\xi$, we get the stated relation
for the second fundamental form of $\D$. This implies also $B^{\D}(\phi X, \phi Y)=
- B^{\D}(X, Y)$ that assures $\mathrm{trace}B^{\D}=0$ (i.e. $\D$ is minimal).

If, in addition, the manifold is K-contact, $\xi$ is Killing, so the induced foliation
$\mathcal{F}_{\xi}$ is Riemannian, which is equivalent to the fact that the orthogonal distribution
$\D$ is totally geodesic.
\end{proof}

The Sasaki condition imposes further restrictions on $B$:

\begin{pr}
Let $(M, \phi, \xi, \eta, g)$ be a Sasaki manifold endowed with an invariant distribution $\mathcal{V}$ which contains $\xi$.
Let $\mathcal{H}$ be the orthogonal complement of $\mathcal{V}$. Then the following
 relations hold:
\begin{equation} \label{bi}
2\left(B^{\mathcal{V}}(U, \phi V) - \phi B^{\mathcal{V}}(U, V) \right)=
\phi(I^{\mathcal{V}}(U, V))-I^{\mathcal{V}}(U, \phi V), \quad \forall U, V \in \Gamma(\mathcal{V}).
\end{equation}
In particular: 
$$
2B^{\mathcal{V}}(U, \xi)+I^{\mathcal{V}}(U, \xi)=0; \quad 
B^{\mathcal{V}}(\phi U, \xi)=\phi  \left(B^{\mathcal{V}}(U, \xi)\right), 
\quad \forall U \in \Gamma (\V).
$$
\end{pr}

\begin{proof}
Note that $\xi\in\Gamma(\mathcal{V})$ implies $\mathcal{H}\subseteq\mathcal{D}$. 
The result now follows from  the definitions and the Sasaki condition: $\nabla_{U}\phi V=
\phi \nabla_{U} V +g(U,V)\xi-\eta(V)U$.

For the second assertion, put $V=\xi$ in formula \eqref{bi} and for the last one, take into
 account the fact that on a Sasaki manifold we have $(\mathcal{L}_{\xi}\phi)X=0$.
\end{proof}

If $\mathcal{V}$ is integrable, we recover the formulas for invariant submanifolds stated in 
\cite[p. 49]{yan}:
\begin{co} 
If $N$ is an invariant submanifold of a Sasaki manifold $M$, then:

$(i)$ $B(X, \xi)=0$

$(ii)$ $B(X, \phi Y)=B(\phi X, Y)=\phi B(X, Y)$
for any vector field X tangent to N (here $B$ denotes 
the second fundamental form of the submanifold). 
\end{co}

\section{Infinitesimal holomorphicity on normal almost contact manifolds}

\subsection{Definitions and first properties}

\begin{de}
Let  $(M, \phi, \xi, \eta)$ be a normal almost contact manifold. A (local) vector field $X$ on $M$ is
called {\it contact - holomorphic} if 
\begin{equation} \label{defhol}
(\mathcal{L}_{X}\phi)Y=\eta \left( [X, \phi Y] \right) \xi, \quad \forall \ Y \in \Gamma(TM).
\end{equation}
A distribution $\mathcal{V}$ on $M$ is called \emph{contact - holomorphic} 
if it admits, around every point, a local frame consisting of contact - holomorphic
 vector fields.
\end{de}

When the context does not impose distinctions, we shall simply write \emph{holomorphic} instead of 
\emph{contact - holomorphic}. 

Holomorphicity of $X$  means collinearity of $(\mathcal{L}_{X}\phi)Y$ 
with $\xi$: the particular form of the coefficient of $\xi$,
generally denoted by $\alpha_{X}(Y)$, results from this 
collinearity condition.

The next result shows the $\phi - invariance$ of the above defined holomorphicity (unlike the usual property
$(\mathcal{L}_{X}\phi)Y=0$):

\begin{lm}
Let $X$ be a holomorphic vector field on  a normal almost contact metric manifold. Then $\phi X$
 is also holomorphic. In particular,  a holomorphic distribution is necessarily invariant.
\end{lm}

\begin{proof}
An explicit formula for the Lie derivative of $\phi$ with respect to $\phi X$ is provided by 
 the following reformulation of the equation \eqref{normalitate}:
$$
(\mathcal{L}_{\phi X}\phi)Y=\phi (\mathcal{L}_{X}\phi)Y - 2d\eta(X,Y)\xi.
$$
Hence, if $X$ holomorphic, then the above equations gives us:
$$
(\mathcal{L}_{\phi X}\phi)Y=-2d\eta(X,Y)\xi.
$$
We now verify that the coefficient of $\xi$ is the same as the one predicted
 by the definition. 
Recall that $\alpha_{X}(Y)=\eta \left( [X, \phi Y] \right)$, so we have
to show that:
$$\alpha_{\phi X}(Y)=\eta \left( [\phi X, \phi Y] \right)=-2d\eta(X,Y).$$
But the normality of $\phi$ assures that
$$
N^{(2)}=0 \Leftrightarrow \eta \left( [\phi X, Y] +[X, \phi Y] \right)=
\phi X \left( \eta(Y) \right)-\phi Y \left( \eta(X) \right).
$$ 
In the above relation we replace $Y$ with $\phi Y$ and we obtain:
$$
\eta \left( [\phi X, \phi Y] -[X, Y]+\eta(Y)[X, \xi]+X(\eta(Y))\xi \right)=
Y \left( \eta(X) \right)-\eta(Y)\xi(\eta(X)),
$$ 
which reduces to
$$
\eta \left( [\phi X, \phi Y]\right) +X(\eta(Y))-Y(\eta(X))-\eta([X, Y])=
-\eta(Y)\left(\xi(\eta(X))-\eta([\xi, X])\right).
$$ 
Finally we use $N^{(4)}:=(\mathcal{L}_{\xi}\eta)X=0$ to derive
$$\eta \left( [\phi X, \phi Y] \right)=-2d\eta(X,Y).$$
\end{proof}

\begin{re}
$(i)$ From the above proof we obtain an alternative expression of the collinearity factor $\alpha_X$:
$$\alpha_{X}(Y)=-\eta \left( [\phi X, Y] \right)+\phi X \left( \eta(Y) \right)-\phi Y \left( \eta(X) \right).$$

$(ii)$ $\alpha_{X}(\xi)=0$ for any holomorphic vector field $X$. This implies that
 $[X, \xi]$ must be collinear with $\xi$ (or, equivalently, $[X^{\mathcal{D}}, \xi]=0$).
In other words, $X$ is projectable with respect to the foliation 
$\mathcal{F}_{\xi}$ locally generated by $\xi$.

$(iii)$ $\alpha_{\xi}(Y)=0$ for any vector field $Y$. Indeed, the normality of $\phi$ implies  
$N^{(3)}:=(\mathcal{L}_{\xi}\phi)Y=0$, so that $\xi$ is holomorphic.

$(iv)$ $X$ is holomorphic if and only if $[X, \xi]$ is collinear with $\xi$ and  
$\left((\mathcal{L}_{X}\phi)Y \right)^{\D}=0$, $\forall Y \in\Gamma(\mathcal{D})$. 
If $M$ is Sasakian, these properties define the \emph{transversally holomorphic fields},
introduced by S. Nishikawa and Ph. Tondeur in \cite {tond}, for manifolds endowed with a
K\"ahler foliation.

\end{re}

\begin{pr} On a normal almost contact manifold, the set $\mathfrak{hol}(M)$ of holomorphic vector fields is a 
Lie subalgebra of $\Gamma(TM)$.
\end{pr}
\begin{proof}
Let $X$ and $X^{\prime}$ be holomorphic vector fields. Then:
\begin{equation*}
\begin{split}
\left(\mathcal{L}_{[X,X^{\prime}]}\phi \right)Y &= \left(\left[\mathcal{L}_{X},\mathcal{L}_{X^{\prime}} \right]\phi \right)Y =
\mathcal{L}_{X}\left(\mathcal{L}_{X^{\prime}} \phi \right)Y-\mathcal{L}_{X^{\prime}}\left(\mathcal{L}_{X} \phi \right)Y\\
&= \left[X,(\mathcal{L}_{X^{\prime}} \phi)Y \right]-(\mathcal{L}_{X^{\prime}} \phi)([X, Y])
- \left[X^{\prime},(\mathcal{L}_{X} \phi)Y \right]+(\mathcal{L}_{X} \phi)([X^{\prime}, Y]).
\end{split}
\end{equation*}
Using the fact that $X$ and $X^{\prime}$ are holomorphic and the remark that $[X, \xi]$ must be
collinear with $\xi$ in this case, we easily obtain that the projection on $\mathcal{D}$ of the
above expression is zero.
Hence $[X, X^{\prime}]$ is also holomorphic.
\end{proof}

On closed Sasakian manifolds with constant transversal scalar curvature,
the structure of $\mathfrak{hol}(M)$ is established in analogy with the K\"ahler case,
cf. \cite {tond}.

\begin{ex}
On $\mathbb{R}^{2n+1}$ with the standard contact metric structure,
take an arbitrary vector field written in an adapted frame as 
$$\displaystyle{X=\alpha \frac{\partial}{\partial z}+\beta^{i}\left(\frac{\partial}{\partial x^{i}}+
y^{i}\frac{\partial}{\partial z}\right)+\gamma^{i}\frac{\partial}{\partial y^{i}}},$$
 where summation is taken with $i=\ov{1,n}$. Note that $ \beta^{i}$ and $\gamma^{i}$ 
are the coefficients of $\displaystyle\frac{\partial}{\partial x^{i}}$ and of 
$\displaystyle\frac{\partial}{\partial y^{i}}$ respectively.
Then $X$ is holomorphic if and only if, for any $i=\ov{1,n}$, $\beta^{i}$ and $\gamma^{i}$ satisfy the Cauchy-Riemann 
equations in the variables $x^{j}$, $y^{j}$ and are constant in $z$: 
$$\displaystyle{\frac{\partial \beta^{i}}{\partial x^{j}}=\frac{\partial \gamma^{i}}{\partial y^{j}}}, \quad 
\displaystyle{\frac{\partial \beta^{i}}{\partial y^{j}}=-\frac{\partial \gamma^{i}}{\partial x^{j}},\quad 
j=\ov{1,n}}, \qquad 
\displaystyle{\frac{\partial \beta^{i}}{\partial z}=\frac{\partial \gamma^{i}}{\partial z}=0}.$$
\end{ex}
The corollary 3.3 below shows that the above description of holomorphic vector fields is not an 
exceptional one. 

\bigskip

As in the complex case (see \cite {moro}, p. 30) we can express the contact-holomorphicity by 
the vanishing of some $\bar\partial$ -operator. In this case 
$\bar\partial: \Gamma(TM) \longrightarrow \mathrm{End}(TM)$ 
satisfies the Leibniz rule and is expressed as follows with respect to Levi-Civita connection:
$$
\bar\partial X(Y)=\frac{1}{2}\phi \left(\nabla_{Y}X + \phi \nabla_{\phi Y}X - \phi(\nabla_{X}\phi)Y \right) 
$$
One can verify that a vector field $X$ is contact-holomorphic if and only if $\bar\partial X(Y)=0$, for all $Y$. 
Equivalently, this means the projectability of $X$ and the vanishing on $X^{\mathcal{D}}$ of a standard 
(transversal) $\bar\partial$ -operator appropriate to $\mathcal{D}$ as $T^{\perp}\mathcal{F}_{\xi}$.
Explicitly: $\bar\partial^{\mathcal{D}}X(Y)=\frac{1}{2}\left(\nabla_{Y}^{\mathcal{D}}X +
\phi \nabla_{\phi Y}^{\mathcal{D}}X - \phi(\nabla_{X}^{\mathcal{D}}\phi)Y \right)$, 
for all $Y\in \Gamma(\mathcal{D})$, where $\nabla^{\mathcal{D}}$ is the adapted connection
in $\mathcal{D}$ in the sense of Tondeur \cite {ton}. Therefore we are dealing with a transversal, projectable notion of holomorphicity 
for vector fields on $M$ regarded as foliated manifold (with the foliation $\mathcal{F}_{\xi}$).

In the Sasaki case, the above formula reduces to:
\begin{equation*}
\begin{split}
\bar\partial X(Y)&=\frac{1}{2}\phi \left(\nabla_{Y}X + \phi \nabla_{\phi Y}X\right), 
\quad \text{for}\,\, Y\in \Gamma(\mathcal{D})\,\, 
\text{and}\\
\bar\partial X(\xi)&=\phi([\xi, X]).
\end{split}
\end{equation*}

\subsection{The holomorphicity condition seen on the cone} 
Recall  that the \emph{cone} $\mathcal{C}(M)$ over an almost contact manifold 
$(M^{2n+1}, \phi, \xi, \eta)$ is $M^{2n+1}\times \mathbb{R}$ with an almost complex 
structure defined by:
$$
J\left(X, f \frac{d}{dt} \right)=(\phi X-f\xi, \eta(X)\frac{d}{dt}).
$$
We point out that the above formula fits the well-known construction of
an almost contact structure on orientable hypersurfaces of almost complex manifolds 
(if we take the standard immersion of $M$ into the cone $\mathcal{C}(M)$ at $t=1$). For details,
 see \cite[Example 4.5.2] {ble}.

\begin{pr}
Let $(M, \phi, \xi, \eta, g)$ be a normal almost contact metric manifold.
As a vector field on the cone over $M$, $\displaystyle (X, f \frac{d}{dt})$ is holomorphic if and only if,
for any $Y \in \Gamma(TM)$, the following relations are satisfied:

\medskip
$(i)$ $(\mathcal{L}_{X}\phi)Y = Y(f)\xi$;

\medskip

$(ii)$ $X(\eta (Y))-\eta \left([X,Y]\right)-\phi Y(f)-\eta(Y) \displaystyle{\frac{df}{dt}}=0$;

\medskip

$(iii)$ $[X, \xi]+\displaystyle{\frac{df}{dt}} \xi=0$;

\medskip

$(iv)$ $\xi(f)=0$.

Hence, if $\displaystyle (X, f \frac{d}{dt})$ is holomorphic on the cone, 
then $X$ is a contact-holomorphic vector field on M.
Moreover, we have the following implications:
\begin{center}
"$(i) \wedge (iii) \Rightarrow (ii)$" and "$(i) \Rightarrow (iv)$". 
\end{center}
\end{pr}

\begin{proof}
One can derive by straightforward computations the following formulas:
\begin{equation*}
\begin{split}
(\mathcal{L}_{(X, f \frac{d}{dt})} J)(Y, 0)&=\left((\mathcal{L}_{X}\phi)Y-Y(f)\xi, (X\eta (Y)-\eta ([X,Y])-\phi Y(f)-\eta(Y) \frac{df}{dt})\frac{d}{dt}\right)\\
(\mathcal{L}_{(X, f \frac{d}{dt})} J)(0, \frac{d}{dt})&=\left(-[X, \xi]-\frac{df}{dt} \xi, \xi(f)\frac{d}{dt}\right)
.
\end{split}
\end{equation*}
As the holomorphicity of $(X, f \frac{d}{dt})$ is equivalent to the vanishing of both expression 
 above, the result follows.

Let us prove the second assertion. 

The implication "$(i) \wedge (iii) \Rightarrow (ii)$" is derived by applying 
 $(i)$ to $\phi Y$ instead of $Y$. We obtain 
$\alpha_{X}(\phi Y)=\phi Y(f)=X\eta (Y)-\eta \left([X,Y]\right)-
\eta(Y) \eta \left([\xi, X] \right)$.
 But from $(iii)$ we have $\eta \left([\xi, X] \right)=\frac{df}{dt}$, 
so the relation $(ii)$ follows.

In order to get "$(i) \Rightarrow (iv)$", put $Y=\xi$ in $(i)$:
 $(\mathcal{L}_{X}\phi)\xi=\xi(f)\xi$. But, as $X$ is holomorphic on $M$, we have already 
noticed that $(\mathcal{L}_{X}\phi)\xi=0$ (i.e. $\alpha_{X}(\xi)=0$), so our implication follows.
\end{proof}

\begin{co}
The contact-holomorphic vector fields on $M$, which come by projection of the 
holomorphic fields on $\mathcal{C}(M)$ form a Lie subalgebra of $\mathfrak{hol}(M)$, 
denoted by $\mathfrak{hol}_{pr}(M)$.
They are contact-holomorphic fields $X$ with two additional properties:

$(a)$ The one-form $\alpha_{X}$ is exact: there exists a function 
$f$ on $M$, such that
$$
Y(f) = \eta \left( [X, \phi Y] \right), \ \forall Y \in \Gamma(TM).
$$

$(b)$ $\eta([X, \xi])$ is (locally) constant (i.e. the factor of collinearity between 
$[X, \xi]$ and $\xi$ is constant).
\end{co}

\begin{proof}
We have seen that, in order to be holomorphic on the cone, a vector field 
must satisfy only $(i)$ and $(iii)$. From condition $(i)$ we obtain $(a)$. 
From $(iii)$, it follows that 
$\displaystyle\frac{df}{dt}=\eta([\xi, X])$, so $f$ is a linear function in $t$: 
$f(p, t)=\eta([\xi, X])t + F(p)$, $p \in M$. In order to verify the equation $(a)$
,
such a function must have the coefficient $\eta([X, \xi])$ (locally) constant, that is $(b)$ holds.

Conversely, if $X$ is a contact-holomorphic vector field on $M$, which satisfies in addition 
$(a)$ and $(b)$, then $\displaystyle \left(X, (\eta([\xi, X])t+f) \frac{d}{dt}\right)$
is holomorphic on $\mathcal{C}(M)$.

In order to see that $\mathfrak{hol}_{pr}(M)$ is a Lie subalgebra, it is enough to note that, on the cone, the holomorphic vector fields 
form a Lie algebra and that the following relation holds:
$$
\left[\left(X, f \frac{d}{dt}\right),\left(X^{\prime}, g \frac{d}{dt}\right)\right]=
\left([X, X^{\prime}], (X(g)-X^{\prime}(f)+f \frac{dg}{dt}-g \frac{df}{dt}) \frac{d}{dt}\right).
$$
\end{proof}

\begin{re}
The subalgebra, $\mathfrak{hol}_{pr}(M)$ contains all vector fields
 along which $\phi$ is invariant: $\mathcal{L}_{X}\phi=0$.
\end{re}

The nature of the constraints $(a)$ and $(b)$ becomes very clear when expressed in local coordinates
 for the case of $\mathbb{R}^{2n+1}$: 

\begin{ex}
On $\mathbb{R}^{2n+1}$ with the standard contact metric structure, let
$\displaystyle X=\alpha \frac{\partial}{\partial z}+\beta^{i} \frac{\partial}{\partial x^{i}}
+\gamma^{i}\frac{\partial}{\partial y^{i}}$ be a holomorphic vector field.

Then $X \in \mathfrak{hol}_{pr}(\mathbb{R}^{2n+1})$ if and only if, in addition, the coefficient 
$\alpha$ takes the form: 
$\alpha=Cz+H(x_{i}, y_{i})$, where $H$ is a harmonic function and $C \in \mathbb{R}$.
\end{ex}

\begin{re}
The relation between contact-holomorphicity on the Sasaki manifolds and holomorphicity 
on its K\"ahler cone can also be obtained 
 using the relations between the Levi-Civita connections on $M$ and $\mathcal{C}(M)$,
 $\nabla$, respectively $\bar{\nabla}$ (for the details, see \cite {gal}). Identifying $X$ on $M$ with $(X,0)$ on the cone, one can prove the following
 formula:
\begin{equation}
(\mathcal{L}_{X} J)Y=(\mathcal{L}_{X} \phi)Y-\left [X(r\eta(Y))+r\eta([X, Y])\right]\partial_{r}
\end{equation}
Indeed, we have the following sequence of identities (where $\Psi:=r\partial_{r}$ is the Euler
field on the cone):
\begin{equation*}
\begin{split}
(\mathcal{L}_{X} J)Y&=\bar{\nabla}_{X}JY-J\bar{\nabla}_{X}Y-\bar{\nabla}_{JY}X+J\bar{\nabla}_{Y}X\\
&=\bar{\nabla}_{X}(\phi Y-\eta(Y) \Psi)-J(\nabla_{X}Y-rg(X, Y) \partial_{r})
-\bar{\nabla}_{\phi Y-\eta(Y) \Psi}X\\
&+J(\nabla_{Y}X-rg(Y, X) \partial_{r})\\
&=\bar{\nabla}_{X} \phi Y-X(\eta(Y)) \Psi-\eta(Y)\bar{\nabla}_{X} \Psi-
J \nabla_{X} Y\\
& -\bar{\nabla}_{\phi Y} X + \bar{\nabla}_{\eta(Y)\Psi} X+J \nabla_{Y} X\\
&=\nabla_{X} \phi Y-rg(X,\phi Y)\partial_{r}-X(\eta(Y)) \Psi-\eta(Y)[X(r)\partial_{r}+r\frac{1}{r}X]\\
&-\phi(\nabla_{X} Y)+\eta(\nabla_{X} Y) \Psi-\nabla_{\phi Y} X+rg(X,\phi Y)\partial_{r}\\
&+\eta(Y)X+\phi(\nabla_{Y} X)-\eta(\nabla_{Y} X) \Psi.
\end{split}
\end{equation*}
This in turn implies  formula (3.2).
\end{re}
\begin{co}
On a normal almost contact metric manifold $(M, \phi, \xi, \eta, g)$ we have:

$(i)$ $a \xi$ is a contact-holomorphic vector field, for any function $a$ defined on $M$
 (so $a \xi \in \mathfrak{hol}(M)$, but not necessarily $a \xi \in \mathfrak{hol}_{pr}(M)$);

$(ii)$ $\displaystyle (\xi, c \frac{d}{dt})$ is a holomorphic vector field on the cone if and only if
 $c$ is a constant.
\end{co}

\begin{proof}
$(i)$ A consequence of normality of $\phi$ (see \cite {ble}) is that $(\mathcal{L}_{\xi} \phi)Y=0$.
 Now, it is an easy task to compute $(\mathcal{L}_{a \xi} \phi)Y=a(\mathcal{L}_{\xi} \phi)Y
-\phi Y(a)\xi$ and to notice that $\alpha_{a \xi}(Y)=-\phi Y(a)\xi$, so the assertion is proved.

$(ii)$ The argument is obvious.
\end{proof}

\subsection{Holomorphicity on Sasakian manifolds}
Recall that on a Riemannian manifold, an arbitrary vector field $V$ induces a derivation $A_{V}$ 
(a tensor field of type $(1, 1)$), defined by: $A_{V}(X):=\nabla_{X}V$. In the complex case,
$A_{V}$ is $J$-linear if and only if $V$ is holomorphic. In our case something similar is happening:
\begin{pr}
On a Sasaki manifold $(M^{2n+1}, \phi, \xi, \eta, g)$ we have:

$(i)$ $V$ is holomorphic if and only if 
$$
\left(A_{V} \circ \phi - \phi \circ A_{V}\right)(X)
 \quad \text{is collinear with} \ \xi, \quad \text{for all}\,\, X \in \Gamma(\mathcal{D})
$$
$\quad$ and also if: $V^{\mathcal{D}}= \phi \nabla_{\xi}V$ (which is equivalent with:
 $[X, \xi]$ collinear with $\xi$).

$(ii)$ If $M^{2n+1}$ is compact and regular and $X$ is a contact-holomorphic vector field on $M^{2n+1}$, 
then $\pi_{*}X$ is holomorphic, where $\pi: M^{2n+1} \longrightarrow M^{2n}$ represents 
the Boothby-Wang fibration. Conversely, the horizontal lift of any holomorphic vector field 
on $M^{2n}$ is a contact-holomorphic vector field on $M^{2n+1}$. 

In particular, the contact distribution on such a Sasaki manifold is holomorphic.
\end{pr}
\begin{proof}
$(i)$ Using the Sasaki condition \eqref{Sasaki} and assuming \eqref{defhol} ($V$ is
 holomorphic), we derive:
$$
\nabla_{\phi X}V=\phi\nabla_{X} V-\eta(X)V +\left(g(V, X)-\eta([V, \phi X])\right)\xi.  
$$
From this, the stated collinearity follows immediately.

Conversely, we can verify that $\eta\left(\nabla_{\phi X}V \right)= g(V, X)-\eta([V, \phi X])$
 and thereafter we can conduct the same calculation backwards to obtain the holomorphicity condition 
 \eqref{defhol}.

$(ii)$ As a direct consequence of the fact that the Boothby-Wang fibration is a Riemannian 
 submersion and satisfies also $\pi_{*}\phi X=J \pi_{*}X$, one get the relation
$$(\mathcal{L}_{\pi_{*}X} J)\pi_{*}Y=\pi_{*}(\mathcal{L}_{X} \phi)Y,$$
 for all projectable vector fields $X, Y$. Note also that (horizontal) contact-ho\-lo\-mor\-phic vector 
 fields on $M^{2n+1}$ 
are, by definition, projectable ($[X^{\mathcal{D}}, \xi]=0$).
The result now follows, as $\xi$ spans $\Ker \pi_{*}$.
\end{proof}

A source of examples of holomorphic vector fields is the following:
\begin{pr}

Let $(M, \phi, \xi, \eta, g)$ be a contact metric manifold. Then any two of the following
conditions imply the third one:

$(i)$ $(\mathcal{L}_{X}g)(Y, Z)=0, \ \forall Y, Z\in\Gamma(\mathcal{D})$,

$(ii)$ $i_{X}d\eta$ is a closed form,

$(iii)$ $X$ is a holomorphic vector field.

Moreover, a vector field $X$ on $M$ is a Killing vector field, which commutes with $\xi$ if and only if $X$ is holomorphic vector
 field which is also strict 
infinitesimal contact transformation (i.e. $\mathcal{L}_{X} \eta = 0$).
\end{pr}
\begin{proof}
The first assertion is a consequence of the following relation:
$$
\left(\mathcal{L}_{X}g\right)(Y, \phi Z)=\left(\mathcal{L}_{X}d\eta\right)(Y, Z)
-g\left(Y, (\mathcal{L}_{X}\phi)Z\right).
$$
For the second assertions we  apply the results obtained by Tanno  
in  \cite[Th. 3.1 and Prop. 3.6]{tano},
 because the holomorphic vector fields which are also strict infinitesimal contact transformation 
are precisely those  for which  $\mathcal{L}_{X} \phi = 0$.
\end{proof}

\begin{re}
The first assertion in the above proposition, can be reformulated in the following terms,
adequate to the foliated structure of the contact metric manifold $M$:
\begin{quote}
{\it a contact - holomorphic vector field with zero transversal divergence is 
a transversal Killing vector field.}
\end{quote}
Clearly, this is a similar result to 
the "if" part of the  Bochner-Yano theorem in the K\"ahler case, cited in 
\cite[p. 93]{kob}. The converse is also true on closed Sasakian manifolds, cf. \cite{tond}.

We recall (in Tondeur's notations, see \cite {ton}) that transversal divergence  of an infinitesimal automorphism of a foliation is 
defined by the relation
$\Theta(X) vol= \mathrm{div}_{B}X \cdot vol$, where $vol$ is a holonomy invariant transversal volume 
($vol= d\eta^{n}$, in our case).
\end{re}

The following analogy with the complex case will be very helpful for local considerations:

\begin{pr}
On a normal almost contact metric manifold $M^{2n+1}$ there always exist (local)
 adapted frames consisting of contact-holomorphic vector fields.
\end{pr}

\begin{proof}
Note first that on the cone over $M$ the vector fields $(\xi, 0)$ 
 and $\displaystyle\left(0, \frac{d}{dt}\right)$ are (real)-holomorphic. Moreover, by construction,
$\displaystyle\left(i \xi, \frac{d}{dt}\right) \in T^{\mathbb{C}}\mathcal{C}(M)$ is a holomorphic vector field 
 on the complexified tangent space to the cone.

 But in our hypothesis, $\mathcal{C}(M)$ is a complex manifold so its tangent bundle is 
holomorphic and then admits local frames of complex holomorphic sections. We can always complete
$\displaystyle \left(i \xi , \frac {d}{dt} \right)$ to such a frame. 

Let 
$\displaystyle\left\{ \left(X_{j}, f_{j}\frac{d}{dt}\right)
-iJ\left(X_{j}, f_{j}\frac{d}{dt}\right) \vert j=\ov{1,n}\right\}$ be such a local completion.

\medskip
We want to prove that $\{\xi, X_{j}^{\mathcal{D}}, \phi X_{j}^{\mathcal{D}} 
\vert j=\ov{1,n}\}$ is an independent family, so it 
represents a local adapted frame for $M$, consisting of contact-holomorphic vector fields.
Observe that $\displaystyle\left(X_{j}, f_{j}\frac{d}{dt}\right)-iJ\left(X_{j}, f_{j}\frac{d}{dt}\right)=
\left(X_{j}-i\phi X_{j}+ if_{j}\xi, (f_{j}-i \eta(X_{j})) \frac{d}{dt} \right)$.

\medskip
Let us now verify that $\left\{X_{j}^{\mathcal{D}}-i\phi X_{j}^{\mathcal{D}} \vert j=\ov{1,n}\right\}$ 
forms a independent family over $\mathbb{C}$, consisting of complex-holomorphic fields. (In the following,
 Einstein convention will be used).
Suppose 
$\lambda^{j} (X_{j}^{\mathcal{D}}-i\phi X_{j}^{\mathcal{D}})=0$. Then we have successively:
$$
\lambda^{j} \left(X_{j}^{\mathcal{D}}-i\phi X_{j}^{\mathcal{D}}, 0 \right)=0,
$$
$$
\lambda^{j} \left(X_{j}-i\phi X_{j}-\eta(X_{j})\xi, 0 \right)=0,
$$
$$
\lambda^{j} \left(X_{j}-i\phi X_{j}+ if_{j}\xi, 0 \right)-
\lambda^{j} ((if_{j}+\eta(X_{j}))\xi, 0 )=0,
$$ 
$$
\lambda^{j} \left(X_{j}-i\phi X_{j}+ if_{j}\xi, (f_{j}-i \eta(X_{j})) \frac{d}{dt} \right)-
\lambda^{j} \left((if_{j}+\eta(X_{j}))\xi, (f_{j}-i \eta(X_{j})) \frac{d}{dt} \right)=0 
$$
and finally
$$
\lambda^{j} \left[\left(X_{j}, f_{j}\frac{d}{dt}\right)-iJ\left(X_{j}, f_{j}\frac{d}{dt}\right)\right]
-\lambda^{j} \left(f_{j}-i\eta(X_{j})\right)(i\xi, \frac{d}{dt})=0.
$$
But this is a linear combination of the vectors that form the complex-holomorphic frame on the cone.
Therefore, $\lambda^{j}=0$ for all $j=\ov{1,n}$.

Now a simple trick will gives us the linear independence over $\mathbb{R}$ of the family 
$\{X_{j}^{\mathcal{D}}, \phi X_{j}^{\mathcal{D}} \mid j=\ov{1,n}\}$.

Suppose that 
$\alpha^{j} X_{j}^{\mathcal{D}} + \beta^{j} \phi X_{j}^{\mathcal{D}}=0$. Then 
$-\beta^{j} X_{j}^{\mathcal{D}} + \alpha^{j} \phi X_{j}^{\mathcal{D}}=0$. 
Together, these relations give 
$\alpha^{j} X_{j}^{\mathcal{D}} + \beta^{j} \phi X_{j}^{\mathcal{D}}
-i(-\beta^{j} X_{j}^{\mathcal{D}} + \alpha^{j} \phi X_{j}^{\mathcal{D}})=0$ which is equivalent to 
$(\alpha^{j}+ i\beta^{j}) (X_{j}^{\mathcal{D}}-i\phi X_{j}^{\mathcal{D}}) =0$, and hence 
$(\alpha^{j}+ i\beta^{j})=0 \Rightarrow \alpha^{j}=\beta^{j}=0$, the relation we wanted to prove.

The argument that $\xi$ is transversal to $\mathcal{D}$ completes the proof.
\end{proof}

A direct computation proves:
\begin{co}
On a normal almost contact manifold, let $\{\xi, E_{i}, \phi E_{i}\}$ be a (local) adapted frame 
 consisting of contact-holomorphic vector fields.
Then a vector field $X=\alpha \xi + \beta^{i} E_{i}+ \gamma^{i} \phi E_{i}$ is holomorphic
if and only if, for all $i=\ov{1,n}$,  $\beta^{i}$ and $\gamma^{i}$ satisfy the \emph{generalized 
Cauchy-Riemann equations}:
$$
E_{j}(\beta^{i})=\phi E_{j}(\gamma^{i}), \quad
E_{j}(\gamma^{i})=-\phi E_{j}(\beta^{i}), \quad  j=\ov{1,n}
$$
and are constant along the flow of $\xi$ (i.e. $\xi(\beta^{i})=\xi(\gamma^{i})=0$).
\end{co}

\subsection{The flow of a contact - holomorphic vector field}

\begin{de}
An map $\psi:(M, \phi, \xi, \eta, g) \longrightarrow (M^{\prime}, \phi^{\prime}, \xi^{\prime}, \eta^{\prime}, g^{\prime})$ 
   between almost contact manifolds is called \emph{contact - holomorphic} if 
$$
d\psi \circ \phi(X)-\phi^{\prime} \circ d\psi(X) \quad \text{is collinear with} \ \xi^{\prime},
 \quad  \forall X \in \Gamma(TM).
$$
As before, the word \emph{contact} in the above notion will be omitted when no confusion is possible. 
\end{de}

\begin{re}
If $\psi$ is holomorphic, then $d\psi(\xi)$ must be collinear with $\xi^{\prime}$. To see this,
put $X=\xi$ in the formula of the above definition.

In particular, the contact-holomorphic maps between
\emph{normal} almost contact manifolds are \emph{transversally holomorphic} as maps between 
foliated manifolds with trans\-ver\-sally holomorphic foliations, according to \cite {drag}
(i.e. $\pi_{\mathcal{D^{\prime}}} \circ d\psi \vert_{\mathcal{D}}$ 
is holomorphic in the usual sense, that is
$(\pi_{\mathcal{D^{\prime}}} \circ d\psi \vert_{\mathcal{D}}) \circ \phi \vert_{\mathcal{D}}=
\phi^{\prime} \vert_{\mathcal{D^{\prime}}} 
\circ (\pi_{\mathcal{D^{\prime}}} \circ d\psi \vert_{\mathcal{D}})$,
where $\pi_{\D}$ stands for the orthogonal projection on the corresponding distribution). 
\end{re}

\begin{pr}
The flow of a contact-holomorphic vector field on a normal almost contact manifold M consists of 
contact-holomorphic transformations on M.
\end{pr}

\begin{proof}
Observe first  that the flow of a holomorphic vector field $\displaystyle\left(X, f \frac{d}{dt} \right)$ on 
 $M \times \mathbb{R}$ decomposes as follows:
$\Psi_{s}=(\psi_{s},\psi_{s}^{t})$, where $\psi_{s}$ can be regarded as the flow of $X$ on $M$ 
 and $\psi_{s}^{t}:M \times \mathbb{R} \longrightarrow \mathbb{R}$, $s\in I_{\epsilon}$
 satisfies the differential equation:
$\displaystyle\frac{d\psi_{s}^{t}}{ds}=f(\psi_{s},\psi_{s}^{t})$.
But we know that if $\displaystyle\left(X, f \frac{d}{dt} \right)$ is holomorphic on the cone 
(which is a complex manifold in this case), 
then its flow $\Psi_{s}$ must be a holomorphic transformation on the cone. We then have successively:
$$
d\Psi_{s} \circ J \left(Y, h \frac{d}{dt} \right)=J \circ d\Psi_{s}\left(Y, h \frac{d}{dt} \right),
$$
$$
 d\Psi_{s}\left(\phi Y-h\xi, \eta(Y)\frac{d}{dt}\right)=
J\left(d\psi_{s}(Y), d\psi_{s}^{t}(Y)+h\frac{\partial\psi_{s}^{t}}{\partial t}\frac{d}{dt}\right),
$$
$$
 \left(d\psi_{s}(\phi Y-h\xi), d\psi_{s}^{t}(\phi Y-h\xi)+\eta(Y)\frac{\partial\psi_{s}^{t}}{\partial t}\frac{d}{dt}\right)=
$$
$$
=\left(\phi(d\psi_{s}(Y))-\left[Y(\psi_{s}^{t})+h \frac{\partial\psi_{s}^{t}}{\partial t} \right]\xi, 
\eta(d\psi_{s}(Y))\frac{d}{dt}\right),
$$
$$
 d\psi_{s}(\phi Y)-\phi(d\psi_{s}(Y))=hd\psi_{s}(\xi)-
\left[Y(\psi_{s}^{t})+h \frac{\partial\psi_{s}^{t}}{\partial t} \right]\xi
$$
and
$$
\phi Y(\psi_{s}^{t})-h\xi(\psi_{s}^{t})+\eta(Y)\frac{\partial\psi_{s}^{t}}{\partial t}=\eta(d\psi_{s}(Y)).
$$
But these two relations must hold also for $Y=0$, that is:
$\displaystyle\psi_{s}(\xi)= \frac{\partial\psi_{s}^{t}}{\partial t} \xi$ and $\xi(\psi_{s}^{t})=0$.
So the above relations reduces to
$$
d\psi_{s}(\phi Y)-\phi(d\psi_{s}(Y))=-Y(\psi_{s}^{t})\xi.
$$
Taking into account that $\xi(\psi_{s}^{t})=0$, the last equation implies,
for $Y=\xi$, that $d\psi_{s}(\xi)$ is collinear with $\xi$.
 
All in all, for the flow of $X$ we have obtained precisely the condition of being 
a contact - holomorphic transformation. 
Moreover we can see 
what means, geometrically, the factor of collinearity with $\xi$.
\end{proof}

\begin{re}
A contact-holomorphic map between Sasakian manifolds is
\emph{transversally harmonic} and an absolute minimum for the energy $E_{T}$
in its foliated homotopy class, according to \cite {drag} (see also \cite {kon}).
\end{re}

\subsection{The $G$-structures viewpoint}
In the end of this section we shall stress out the connection between a certain $G$-structure of 
almost contact manifolds and the contact-holomorphicity, which have been discussed until now
(for general definitions, see \cite {kob}).

The existence of an almost contact (metric) structure on a manifold $M^{2n+1}$ is equivalent with the 
existence of a 
$(U(n) \times 1)$-structure which clearly is not integrable (even when $\phi$ is normal).
The normality of $\phi$ reflects in the integrability of other $G$-structure of $M^{2n+1}$, namely the 
$H^{1,n}$-structure, called also \emph{transversal holomorphic structure} 
 (for notations and details, see \cite {gom}). 
The infinitesimal automorphisms of the $H^{1,n}$-structure are precisely the 
contact - holomorphic vector fields that we have dealt with, so far. In a system of (local) 
 distinguished coordinates $(u, z^{j}, \ov{z}^{j})$, these vector fields take the form
\begin{equation*}
\begin{split}
&X = a(u, z^{j}, \ov{z}^{j})\frac{\partial}{\partial u}+
b_{k}(u, z^{j}, \ov{z}^{j})\frac{\partial}{\partial z^{k}}+
\ov{b}_{k}(u, z^{j}, \ov{z}^{j})\frac{\partial}{\partial \ov{z}^{k}},\,\, \text{where}\\
&\frac{\partial b_{k}}{\partial \ov{z}^{j}}=0 \quad \text{and}\quad \frac{\partial b_{k}}{\partial u}=0.
\end{split}
\end{equation*}

If, in addition, $M^{2n+1}$ is contact, passing from these coordinates to Darboux coordinates 
 will not respect the $H^{1,n}$-structure,
 so the distinguished coordinates and above local expression for $X$ will be not at all suited for the 
study of strict contact geometric properties of $M^{2n+1}$.  

\section{Complex holomorphicity on normal almost contact manifolds}

In this section we stress out the notion of holomorphic vector field
 in the complex context. If $(M, \phi, \xi, \eta, g)$ is a normal almost contact metric 
manifold, then the complexified tangent bundle admits a natural split:
$$
T^{\mathbb{C}}M=T^{0}M \oplus T^{(1,0)}M \oplus T^{(0,1)}M,
$$
 where $T^{(1,0)}M=\{X-i \phi X \mid X \in \Gamma(\mathcal{D})\}$, $T^{(0,1)}M=\ov{T^{(1,0)}M}$
 and $T^{0}M=Sp_{\mathbb{C}}\{\xi\}$ are the eigenspaces of $\phi$ corresponding to
 the eigenvalues $i, -i$ and 0.

\begin{de}
On an almost contact manifold $(M, \phi, \xi, \eta)$, a smooth function 
$f: M \longrightarrow \mathbb{C}$ will be called \emph{holomorphic} if $df \circ \phi=i \cdot df$.
\end{de}

\begin{pr}
Let $f:M \longrightarrow \mathbb{C}$ be a smooth function on a normal almost contact manifold
 M. Then the following statements are equivalent:

$(i)$ $f$ is holomorphic,

$(ii)$ $Z(f)=0$, for all $Z\in T^{0}M \oplus T^{(0,1)}M$,

$(iii)$ $df \in \Lambda^{(1, 0)}_{B}M$, where $\Lambda^{(1, 0)}_{B}M$ comes from the natural
splitting of the complexification of the sheaf of basic one forms on M:
$\Lambda^{1}_{B} \otimes \mathbb{C}=\Lambda^{(1, 0)}_{B} \oplus \Lambda^{(0, 1)}_{B}$,
cf. \cite {bgn}.

In addition, if $\psi: M \longrightarrow M$ is a holomorphic map, then $f \circ \psi$ 
is a holomorphic function on M.
\end{pr}

\begin{proof}
In order to prove "$(i) \Leftrightarrow (ii)$", we have simply to remark that $df(\xi)=0$ 
(so $\xi(f)=0$) and then the rest of the proof will be similar to the complex case:

$df (\phi X) =i df(X) \Leftrightarrow i df (X + i\phi X)=0 \Leftrightarrow (X + i\phi X)(f)=0,
\forall X \in \Gamma(TM)$.

In the proof of "$(i) \Leftrightarrow (iii)$" it suffices to verify that $df$ 
is a basic 1-form. We have already seen that $df(\xi)=0$. It remains to compute: 
$$
(\mathcal{L}_{\xi}df)(X)=\xi(df(X))-df([\xi, X])=\xi(X(f))-[\xi, X](f)=X(\xi(f))=0.
$$

For the last assertion, we have to do a simple verification:

$d(f \circ \psi)(\phi X)=df(d\psi(\phi X))=df(\phi(d\psi(X))+a\xi)=df(\phi(d\psi(X))=idf(d\psi(X))$.
\end{proof}

\begin{de}
On a normal almost contact metric manifold $M$, $Z\in T^{0}M \oplus T^{(1,0)}M$ will be called
 \emph{complex - holomorphic} if $Z(f)$ is holomorphic for any (local) holomorphic function $f$ on $M$.
\end{de}

\begin{pr}
$Z=a \xi+X-i\phi X \in T^{0}M \oplus T^{(1,0)}M$ is complex-holomorphic if and only if X is holomorphic 
$($in the expression of Z, a is a complex valued function and $X \in \Gamma(\mathcal{D}))$.
\end{pr}

\begin{proof}
Let $Z=a \xi+X-i\phi X$ be a complex-holomorphic vector field and $f$ a holomorphic function
 on $M$. We have seen that $(X+i\phi X)(f)=0$, so $Z(f)=(X-i\phi X)(f)=2X(f)$ must be a
holomorphic function. This means also that: $(Y+i\phi Y)(X(f))=0, \forall Y \in TM$.

From all this we can deduce that: $[Y+i\phi Y, X](f)=0$ (for an arbitrary holomorphic function $f$),
 which in turn implies:
 $[Y+i\phi Y, X]\in T^{0}M \oplus T^{(0,1)}M$. 

But, for any
 $W=a \xi+Y+i\phi Y\in T^{0}M \oplus T^{(0,1)}M$, we have: $\mathrm{Im}(W)^{\mathcal{D}}
=\phi(\mathrm{Re}(W)^{\mathcal{D}})$.
 In our case, $\mathrm{Im}\left([Y+i\phi Y, X]\right)=[\phi Y,X]$ and
$\mathrm{Re}\left([Y+i\phi Y, X]\right)=[Y,X]$.
So we must have: 
$$
[\phi Y,X]^{\mathcal{D}}=\phi([Y,X]^{\mathcal{D}}) \Leftrightarrow
((\mathcal{L}_X \phi)Y)^{\mathcal{D}}=0
$$
and this means that $X$ is holomorphic.

\medskip
Conversely, let $X$ be a holomorphic vector field and $f$ a holomorphic function. We have to show
 that $Z(f)=(a \xi+X-i\phi X)(f)$ is a holomorphic function too. But
$Z(f)=(X-i\phi X)(f)=2X(f)$, because $\xi(f)=(X+i\phi X)(f)=0$, $f$ being holomorphic.
According to Prop. 4.1, $X(f)$ is holomorphic if and only if $(b \xi+Y+i\phi Y)(X(f))=0$,
for any $b$ complex valued function and $Y \in \Gamma(\mathcal{D})$. In turn, this is equivalent 
to $[b \xi+Y+i\phi Y, X](f)=0$ which is assured by $[b \xi+Y+i\phi Y, X]\in T^{0}M \oplus T^{(0,1)}M$
(due to the holomorphicity of $X$).

\end{proof}

Analogous as in the complex case, we have also:

\begin{pr}
On a normal almost contact metric manifold, $T^{0}M \oplus T^{(1,0)}M$ and $T^{0}M \oplus T^{(0,1)}M$ are
 integrable subbundles of $T^{\mathbb{C}}M$, invariant along the flow of a holomorphic
 vector field X.
\end{pr}

\begin{proof}
We have to prove that $[a \xi+X-i\phi X, b \xi+Y-i\phi Y] \in T^{0}M \oplus T^{(1,0)}M$.

A well known result of Ianu\c s, \cite {sian}, tells us that, in this case, $T^{(1,0)}M$ is involutive.
 So it remains to prove that $[X-i\phi X, b \xi] \in T^{0}M \oplus T^{(1,0)}M$.

Taking into account that $\mathcal{L}_{\xi}\phi=0$ (i.e. $[\xi, \phi X]=\phi [\xi, X], \forall X$), we have:
\begin{equation*}
\begin{split}
[X-i\phi X, b \xi] & = (X-i\phi X)(b)\xi+ b \left( [X, \xi]-i[\phi X, \xi] \right) \\
& = (X-i\phi X)(b) \xi +b \left([X, \xi]-i \phi([X, \xi]) \right) \\
&\in T^{0}M \oplus T^{(1,0)}M.
\end{split}
\end{equation*}

As usual $\psi_{s}$ denote the flow of $X$. We have:
\begin{equation*}
\begin{split}
d\psi_{s}(a \xi+X-i\phi X)&=a d\psi_{s}(\xi)+d\psi_{s}(X)-id\psi_{s}(\phi X)=\\
&=a b\xi+d\psi_{s}(X)-i\left(\phi (d\psi_{s}X)+b^{\prime} \xi \right)=\\
&=(a b-ib^{\prime})\xi+d\psi_{s}(X)-i\phi (d\psi_{s}X)\\
&\in T^{0}M \oplus T^{(1,0)}M.
\end{split}
\end{equation*}
\end{proof}

\begin{re}
Note that in Proposition 3.5. we have proved that $T^{0}M \oplus T^{(1,0)}M$ admits, locally,
frames of holomorphic sections.
\end{re}

The proof of the following proposition is an easy computation and we shall omit it:

\begin{pr}
On a Sasaki manifold we always have:

\medskip
$(i)$ $\nabla_{\ov{W}}Z \in T^{0}M \oplus T^{(1,0)}M, \quad \forall W, Z \in T^{(1,0)}M.$

\medskip
$(ii)$ $\nabla_{a \xi}Z \in T^{(1,0)}M, \quad \forall Z \in T^{(1,0)}M$.

\medskip
$(iii)$ $\nabla_{\ov{W}}a \xi \in T^{0}M \oplus T^{(0,1)}M, \quad \forall W \in T^{(1,0)}M.$

\medskip
In addition, $Z \in T^{(1,0)}M$ is a complex holomorphic field if and only if:
$$
\nabla_{\ov{W}}Z \in T^{0}M, \ \forall W \in T^{(1,0)}M \quad \text{and} \quad \nabla_{\xi}Z=-iZ.
$$
\end{pr}

\begin{re}
The contact (complex) holomorphicity, which we deal with, is more general
than the one introduced by Tanaka in \cite {tan}. One can verify that a contact complex-holomorphic
field from $T^{(1,0)}M$ is holomorphic also in Tanaka's sense if, in addition, it preserves
the contact distribution, 
or, equivalently, if $\phi$ is invariant along its flow (i.e. $\mathcal{L}_{X}\phi=0$).
This is a rather strong restriction (generally not satisfied in our context).
\end{re}

\section{Holomorphic foliations on a Sasaki manifold}

Again by analogy with the K\"ahler case (treated in \cite {sve}), in the following 
we shall stress out some properties of the 
holomorphic distributions. For the sake of completeness we recall the notion of 
\emph{mixed sectional curvature} of a Riemannian manifold $M$ endowed with two complementary distributions $\V$ and $\Hh$:
$$
s_{mix}=\sum_{i, \alpha} K^{M}(e_{i} \wedge f_{\alpha})
$$
where $\{e_{i}\}, \{f_{\alpha}\}$ are local orthonormal frames for $\mathcal{V}$ and 
$\mathcal{H}$.

\begin{pr}
On a Sasaki manifold $(M^{2n+1}, \phi, \xi, \eta, g)$, an invariant holomorphic distribution
 $\mathcal{V}$ of dimension $2p+1$ has the following properties $($as usual, $\mathcal{H}=\mathcal{V}^{\bot})$:

\medskip
$(i)$ $\mathcal{V}(\nabla_{\phi Z}X+\nabla_{Z}\phi X)=0, \quad \forall Z \in\Gamma(TM), X \in\Gamma(\mathcal{H}).$

\medskip
$(ii)$ $\phi B^{\mathcal{H}}(X, Y) + g(X, Y)\xi =\frac{1}{2}I^{\mathcal{H}}(X, \phi Y), \quad \forall X, Y \in \Gamma(\mathcal{H}).$

\medskip
$(iii)$ $\vert B^{\mathcal{H}}\vert^{2}+2(n-p)=\frac{1}{4}\vert I^{\mathcal{H}}\vert^{2}$.

\medskip
$(iv)$ $\mathrm{trace}B^{\mathcal{V}}=0$ $(\mathcal{V}$ is a minimal distribution$)$.
\end{pr}

\begin{proof}
$(i)$ Because $M$ is Sasakian, we have:
$(\nabla_{V}\phi)Z=g(V,Z)\xi-\eta(Z)V$. So, for any section $X$ of $\mathcal{H}$ and $V$ of 
$\mathcal{V}$, the following relation holds: $g \left(X, (\nabla_{V}\phi)Z \right)=0$,
also because $\xi \in \Gamma(\mathcal{V})$, by Prop. 2.1. 
Taking this into account, together with the holomorphicity hypothesis, we derive the relation $(i)$ using
Prop. 2.2.

$(ii)$ Using $(i)$, we have:

\medskip
\noindent
$g \left(\frac{1}{2}I^{\mathcal{H}}(X, \phi Y), V \right)=g \left(\frac{1}{2}(\nabla_{X}\phi Y-
\nabla_{\phi Y}X), V \right)=g \left(\frac{1}{2}(\nabla_{X}\phi Y+\nabla_{Y}\phi X), V \right)=$

\medskip
\noindent
$=\frac{1}{2}g \left(\phi \nabla_{X} Y+g(X,Y)\xi-\eta(Y)X+\phi \nabla_{Y} X+g(Y,X)\xi-\eta(X)Y, V \right)=$

\medskip
\noindent
$=\frac{1}{2}g \left(\phi (\nabla_{X} Y+\nabla_{Y} X)+2g(X,Y)\xi-\eta(Y)X-\eta(X)Y, V \right)=$

\medskip
\noindent
$=g \left(\phi B^{\mathcal{H}}(X, Y), V \right)+g(X,Y)g(\xi, V).$ 

The last equality completes the proof because all the terms in the relation 
$(ii)$ are sections of $\V$ and $V \in \mathcal{V}$ was arbitrary.

$(iii)$ This formula involving the Hilbert-Schmidt norms of $B^{\mathcal{H}}$ and 
$I^{\mathcal{H}}$ is a straight consequence of $(ii)$ if we point out that:

$\eta (B^{\mathcal{H}}(X, Y))=g \left(B^{\mathcal{H}}(X, Y), \xi \right)=- \frac{1}{2}(\mathcal{L}_{\xi} g)(X, Y)=0$, 
because $\xi$ is a Killing
 vector field in the Sasakian context.

This assures that $\Vert \phi B^{\mathcal{H}}(X, Y) \Vert = \Vert B^{\mathcal{H}}(X, Y) \Vert$.

In order to compute $\vert I^{\mathcal{H}}\vert^{2}$, it is worth to notice that 
$\xi \in \Gamma(\mathcal{V})$ implies $\mathcal{H} \subseteq \mathcal{D}$. 
So, for a local frame of 
$\mathcal{H}$ of the type $\{e_{i}, \phi e_i \}$, we shall have: $\phi^{2}e_{i}=-e_{i}$.

$(iv)$ The relation \eqref {bi} can be rewritten as follows:
$$
2\left(B^{\mathcal{V}}(U, \phi V) - \phi B^{\mathcal{V}}(U, V) \right)=
-\left[(\mathcal{L}_{U}\phi) V\right]^{\mathcal{H}}, \quad \forall \ U, V \in \Gamma(\mathcal{V}).
$$

For a (contact-)holomorphic field $U$, we get: $B^{\mathcal{V}}(U, \phi V) = 
\phi B^{\mathcal{V}}(U, V)$, which implies immediately $B^{\V}(U, \phi V) = B^{\V}(\phi U, V)$. 

Using also that $[U, \xi]$ is collinear with $\xi$ when $U$ is holomorphic (so $I^{\V}(U, \xi)=0$),
again from Prop. 2.3 we obtain:
$$
B^{\V}(U, V)+B^{\mathcal{V}}(\phi U, \phi V) = 0, \quad \forall \ U, V \in \mathfrak{hol}(M).
$$

Therefore, in a local frame of holomorphic vector fields, we will have:
$$
\mathrm{trace}B^{\mathcal{V}}=\nabla_{ \xi} \xi+\sum_{i}\mathcal{H}\left[\nabla_{e_{i}} e_{i}+
\nabla_{\phi e_{i}} \phi e_{i} \right]=0.
$$
\end{proof}

\begin{pr}
Under the same hypothesis as above, Walczak formula $($see \cite {wal}$)$ simplifies to:
\begin{equation}\label{wal}
\mathrm{div}^{\mathcal{V}}\mathrm{trace} B^{\mathcal{H}}+2(n-p)+
\frac{1}{4}\vert I^{\mathcal{V}}\vert^{2}= s_{mix}+\vert B^{\mathcal{V}}\vert^{2}
\end{equation}
\end{pr}

\begin{proof}
Recall that, for an arbitrary Riemannian manifold $(M, g)$ with two orthogonal complementary 
distributions $\mathcal{V}$ and $\mathcal{H}$, Walczak formula asserts:
\begin{equation*}
\begin{split}
& \mathrm{div}^{\mathcal{V}}\mathrm{trace} B^{\mathcal{H}}+
\mathrm{div}^{\mathcal{H}}\mathrm{trace} B^{\mathcal{V}}+
\frac{1}{4}\vert I^{\mathcal{H}}\vert^{2}+
\frac{1}{4}\vert I^{\mathcal{V}}\vert^{2} \\
 &= s_{mix}+\vert B^{\mathcal{H}}\vert^{2}+
\vert B^{\mathcal{V}}\vert^{2},
\end{split}
\end{equation*}
Now, applying $(iii)$ and $(iv)$ from Proposition 5.1, the result follows. 
\end{proof}

\begin{re}
When $\mathcal{V}$ is integrable, the equation \eqref{wal} reduces to:
\begin{equation}\label{wa}
\mathrm{div}^{\mathcal{V}}\mathrm{trace}B^{\mathcal{H}}+2(n-p)=
s_{mix}+\vert B^{\mathcal{V}}\vert^{2}
\end{equation}
\end{re}

Integrating \eqref{wa} along any compact leaf, we get the following:
\begin{te}
{\bf Bochner-type result}

Let $(M^{2n+1}, \phi, \xi, \eta, g)$ be a Sasaki manifold with a $(2p+1)$-dimensional holomorphic
 foliation such that $s_{mix}\geq 2(n-p)$. Then $s_{mix}= 2(n-p)$ along every compact leaf and every compact leaf
 is a totally geodesic submanifold of $M$. In particular, if $s_{mix}> 2(n-p)$, then
 $\mathcal{V}$ cannot have compact leaves.
\end{te}

\begin{co}
Let $(M^{2n+1}, \phi, \xi, \eta, g)$ be a compact Sasaki manifold with the sectional curvature
 $k \geq 2m$ ($m < n$). Then every $(\phi, J)$-holomorphic submersion from $M$ into any Hermitian
 manifold $N^{2m}$ has totally geodesic fibers.
\end{co}

Other results as Prop. 3.8. and Prop. 3.9. in \cite {sve}, dealing with holomorphic conformal
 foliations, can be also restated, now in a obvious way, for the Sasakian case.
 
It is worth to notice that the $(\phi, J)$-holomorphic submersions on Sasaki manifolds
into a K\"ahler manifold are in fact a special class of pseudo-harmonic morphisms, with very nice
geometric properties, cf. \cite {aab}. 

\begin{pr}
Let $(M^{2m+1}, \phi, \xi, \eta, g)$ be a Sasaki manifold.
Then every $(\phi, J)$-holomorphic submersion $\psi$, from $M$ onto a K\"ahler manifold
$(N^{2n}, J, g_{N})$ is a pseudo-horizontally homothetic (PHH) harmonic morphism.

In particular, it has minimal fibers and the inverse images of complex submanifolds
in $N$ are invariant, so minimal, submanifolds of $M$.
If in addition $m=n$, then the horizontal distribution (of the submersion $\psi$),
$\mathcal{H}$, coincides with the contact distribution on
$M$ (in particular $\mathcal{H}$ cannot be integrable).
\end{pr}

\begin{proof}
The harmonicity of such submersions has been remarked already in \cite {ian}. Then we have to verify
the PHWC condition (Pseudo Horizontal Weak Conformality) and the PHH one.

The first condition
simply means that the induced almost complex structure on the horizontal bundle (defined by
 $J_{\mathcal{H}}=d\psi^{-1} \circ J \circ d\psi$) is compatible with the metric $g$.
That is indeed the case, because $\mathcal{H} \subset \mathcal{D}$ 
(due to $\xi \in \Gamma(\mathcal{V})$) and $J_{\mathcal{H}}$ coincides 
with $\phi$ restricted to $\mathcal{H}$ (due to the $(\phi, J)$-holomorphicity of $\psi$). 

The second (PHH) condition means that $J_{\mathcal{H}}$ is parallel in horizontal directions
with respect to $\nabla^{\mathcal{H}}$, so it satisfies a \emph{partial K\"ahler} condition.
To see this we have to particularize the formula (1.2) for $X, Y \in \Gamma(\mathcal{H})
 \subset \Gamma(\mathcal{D})$ and to take the $\mathcal{H}$-part of both sides of the relation.
 
\end{proof}

{\bf Acknowledgements}.
This work has been done during the visit at \emph{Ecole Polytechnique F\'ed\'erale
 de Lausanne} in the frame of {\bf SCOPES Programme}. The authors are grateful to
 Tudor Ra\c tiu for his hospitality and for the fruitful discussions that they had.

The second named author thanks  Liviu Ornea for constant encouragements and constructive discussions.


\begin{thebibliography}{100}

\bibitem {aab} M.A. Aprodu, M. Aprodu, V. Brinzanescu, \emph{A class of harmonic submersions and minimal submanifolds}, 
Int. J. of Math. {\bf 11}(9) (2000) 1177-1191.
 
\bibitem {drag} E. Barletta, S. Dragomir, \emph{On transversally holomorphic maps of K\"ahlerian foliations}, 
Acta Applicandae Mathematicae, {\bf 54}(1998), 121-134.

\bibitem {ble} D.E. Blair, \emph{Riemannian Geometry of Contact and Symplectic Manifolds}, Birkhauser Boston, 
Progress in Mathematics, vol.203, 2002.

\bibitem {gal} Ch. Boyer, K. Galicki, \emph{3-Sasakian Manifolds}, in \emph{Surveys in differential geometry: 
essays on Einstein manifolds}, in Surv. Differ. Geom., VI, Int. Press, Boston, MA, 1999, 123-184.

\bibitem {bgn} Ch. Boyer, K. Galicki, M. Nakamaye, \emph{On the geometry of Sasakian-Einstein 5-manifolds}, 
Math. Ann. {\bf 325}(3) (2003), 485-524

\bibitem {gom} X. Gomez-Mont, \emph{Transversal holomorphic structures}, J. Diff. Geom. {\bf 15} (1980), 161-185.

\bibitem{sian} S. Ianu\c s, \emph{Sulle varieta di Cauchy-Rieman}, Rend. dell'Accademia di Scienze Fisiche e Matematiche, 
Napoli, {\bf XXXIX}(1972), 191-195.

\bibitem {ian} S. Ianu\c s, A.M. Pastore, \emph{Harmonic maps on contact metric manifolds}, Ann. Math. Blaise Pascal {\bf 2} (1995),
 43-53.

\bibitem {kob} S. Kobayashi, \emph{Transformation Groups in Differential Geometry}, Springer-Verlag (1972).

\bibitem {kon} J. Konderak, R. Wolak, \emph{Transversally harmonic maps between manifolds with Riemannian foliations}, 
Q.J. Math. {\bf 54}(2003), 335-354.

\bibitem {moro} A. Moroianu, \emph{Lectures on K\"ahler Geometry}, arXiv:math.DG/0402223.

\bibitem {tond} S. Nishikawa, Ph. Tondeur, \emph{Transversal infinitesimal automorphisms for harmonic K\"ahler foliations},
 T\^ohoku Math. J. {\bf 40}(1988), 599-611.

\bibitem {ol} Z. Olszak, \emph{On contact metric manifolds}, Tohoku Math. J., {\bf 31} (1979), 247-253. 

\bibitem {sve} M. Svensson, \emph{Holomorphic foliations, harmonic morphisms and the Walczak formula}, 
J. London Math. Soc.(2) {\bf 68}(3) (2003), 781-794.

\bibitem {tan} N. Tanaka, \emph{A Differential Geometric Study on Strongly Pseudo-Convex Manifolds}, Lectures in Mathematics, 
KYOTO University, 1975.

\bibitem {tano} S. Tanno, \emph{Some transformations on manifolds with almost contact and contact metric structures I}, 
Tohoku Math. J. {\bf 15} (1963), 140-147.

\bibitem {ton} Ph. Tondeur, \emph{Foliations on Riemannian manifolds},Universitext, Springer Verlag, New York, 1988.

\bibitem {yan} K. Yano, M. Kon, \emph{CR Submanifolds of Kaehlerian and Sasakian Manifolds} Birkhauser Boston, 
Progress in Mathematics, vol. {\bf 30}, 1983.
 
\bibitem {wal} P.G. Walczak, \emph{An integral formula for a Riemannian manifold with two orthogonal complementary distributions}, 
Colloquium Mathematicum {\bf 58}(2) (1990), 243-252.

\end{thebibliography}
\end{document}